\documentclass{article}%
\usepackage{amsfonts}
\usepackage{amsmath}
\usepackage{hyperref}
\usepackage{amssymb}
\usepackage{graphicx}
\usepackage{theorem}
\usepackage{accents}%
\newtheorem{theorem}{Theorem}

\newtheorem{conjecture}[theorem]{Conjecture}

\newtheorem{remark}{Remark}

\begin{document}

\title{Bernoulli Partitions}
\author{Thomas Curtright\medskip\\Department of Physics, University of Miami, Coral Gables, FL 33124\\{\footnotesize curtright@miami.edu}}
\date{}
\maketitle

\begin{abstract}
Scale invariant scattering suggests that all Bernoulli numbers $B_{2n}$ can be
naturally partitioned, i.e., written as particular finite sums of same-signed,
monotonic, rational numbers. \ Some properties of these rational numbers are
discussed here, especially in the limit of large $n$.

\end{abstract}

\begin{center}
\textit{In tribute to Luca Mezincescu (1946-2025)}
\end{center}

\noindent\hrulefill\bigskip

In this note I extend some earlier observations \cite{TLCSigma} about
Bernoulli numbers as obtained in the context of computing the scattering
amplitude for a scale invariant potential \cite{TLCJMP,TLCBJP}.

In two spatial dimensions (2D), when computed using quantum mechanics,
non-relativistic scattering by a repulsive inverse square potential
$V=\kappa/r^{2}$ ($\kappa>0$) results in a simple form for the integrated
cross section, $\sigma=\int_{0}^{2\pi}\left(  \frac{d\sigma}{d\theta}\right)
d\theta$. \ For a mono-energetic beam of mass $m$ particles the 2D result is
$\sigma=2\pi^{2}m\kappa/\left(  \hbar^{2}k\right)  $ where the incident energy
is $E=\hbar^{2}k^{2}/\left(  2m\right)  $. \ This result follows from a
straightforward application of 2D phase-shift analysis for the potential
$V=\kappa/r^{2}$ upon realizing a remarkable identity involving the sinc
function, $\operatorname{sinc}\left(  z\right)  \equiv\sin\left(  z\right)
/z$. \ 

A succinct form of the identity in question is \cite{Illinois}
\begin{equation}
1=\frac{\sin\left(  \pi x\right)  }{\pi x}+2\sum_{l=1}^{\infty}\frac{\left(
-1\right)  ^{l}\sin\left(  \pi\sqrt{l^{2}+x^{2}}\right)  }{\pi\sqrt
{l^{2}+x^{2}}} \label{RID1}%
\end{equation}
All higher powers of $x$ cancel when terms on the RHS are expanded as series
in $x^{2}$, as a consequence of familiar $\zeta\left(  2n\right)  $ exact
values for integer $n>0$, as shown in \cite{TLCJMP}. \ Upon expressing the
sinc function in terms of spherical Bessel functions, and then using series
representations for Bessel functions \cite{AandS} in terms of Bernoulli
numbers $B_{n}$, the identity (\ref{RID1}) leads directly to a set of finite
sum identities involving those numbers. \ Namely,%
\begin{equation}
1=\left(  -1\right)  ^{n+1}\left(  4n+2\right)  \sum_{k=0}^{n}\frac{\left(
2n\right)  !}{n!k!\left(  n-k\right)  !}\left(  \frac{B_{n+k+1}}%
{n+k+1}\right)  \text{ \ \ for integer }n\geq1\text{.} \label{RID2}%
\end{equation}
This finite sum involves
\href{https://en.wikipedia.org/wiki/Trinomial_expansion}{trinomial
coefficients} as well as
\href{https://en.wikipedia.org/wiki/Bernoulli_number}{divided Bernoulli
numbers}, $\beta_{m}\equiv B_{m}/m$. \ It is not difficult to check the
validity of (\ref{RID2}) using various expressions of the Bernoulli numbers as
finite, \emph{alternating} sums \cite{Reciprocity}. \ Recall the usual phases
are given by $B_{2n}=\left(  -1\right)  ^{n+1}\left\vert B_{2n}\right\vert $
as well as $B_{2n+1}=0$ for $n=1,2,\cdots$. \ Also recall the well-known
relation
\begin{equation}
\left\vert B_{2n}\right\vert =\frac{2\left(  2n\right)  !}{\left(
2\pi\right)  ^{2n}}~\zeta\left(  2n\right)  \text{ \ \ for \ \ }%
n=1,2,\cdots\label{well-known}%
\end{equation}
Thus, given an identity for $\left\vert B_{2n}\right\vert $ a similar identity
for $\zeta\left(  2n\right)  $ follows immediately.

As quipped previously \cite{TLCSigma}, if encountered as graffiti on the
stones of a bridge (say, in Ireland) either (\ref{RID1}) or its companion
identity (\ref{RID2}) might cause nothing more than a raised eyebrow in
passing. \ \newpage

\noindent Perhaps justifiably so. \ However, upon inverting the linear
relations in (\ref{RID2}) to obtain expressions for each individual Bernoulli
number, the results are more striking: \ The Bernoulli numbers $B_{2n}$ for
$n\geq2$ are \emph{naturally partitioned}. \ That is to say, each $B_{2n}$ is
given by an interesting sum of $n-1$ monotonic, same-signed, rational numbers.
\ Unlike many well-known representations of $B_{2n}$, here the terms in the
finite sums do \emph{not} alternate in sign \cite{Manifesto}.

For example, the first ten partitions are given by \
\[
\left(
\begin{array}
[c]{c}%
\left\vert B_{2}\right\vert \smallskip\\
\left\vert B_{4}\right\vert \smallskip\\
\left\vert B_{6}\right\vert \smallskip\\
\left\vert B_{8}\right\vert \smallskip\\
\left\vert B_{10}\right\vert \smallskip\\
\left\vert B_{12}\right\vert \smallskip\\
\left\vert B_{14}\right\vert \smallskip\\
\left\vert B_{16}\right\vert \smallskip\\
\left\vert B_{18}\right\vert \smallskip\\
\left\vert B_{20}\right\vert \smallskip
\end{array}
\right)  =\left(
\begin{array}
[c]{c}%
\frac{1}{6}\smallskip\\
\frac{1}{30}\smallskip\\
\frac{1}{42}\smallskip\\
\frac{1}{30}\smallskip\\
\frac{5}{66}\smallskip\\
\frac{691}{2730}\smallskip\\
\frac{7}{6}\smallskip\\
\frac{3617}{510}\smallskip\\
\frac{43\,867}{798}\smallskip\\
\frac{174\,611}{330}\smallskip
\end{array}
\right)  =\left(
\begin{array}
[c]{l}%
\frac{1}{6}\smallskip\\
\frac{1}{30}\smallskip\\
\frac{1}{60}+\frac{1}{140}\smallskip\\
\frac{1}{45}+\frac{1}{105}+\frac{1}{630}\smallskip\\
\frac{1}{20}+\frac{3}{140}+\frac{1}{252}+\frac{1}{2772}\smallskip\\
\frac{1}{6}+\frac{1}{14}+\frac{17}{1260}+\frac{1}{693}+\frac{1}{12\,012}%
\smallskip\\
\frac{691}{900}+\frac{691}{2100}+\frac{59}{945}+\frac{41}{5940}+\frac
{5}{10\,296}+\frac{1}{51\,480}\smallskip\\
\frac{14}{3}+\ 2\ +\frac{359}{945}+\frac{8}{189}+\frac{4}{1287}+\frac{1}%
{6435}+\frac{1}{218\,790}\smallskip\\
\frac{3617}{100}+\frac{10\,851}{700}+\frac{1237}{420}+\frac{217}{660}%
+\frac{293}{12\,012}+\frac{1}{780}+\frac{7}{145\,860}+\frac{1}{923\,780}%
\smallskip\\
\frac{43\,867}{126}+\frac{43\,867}{294}+\frac{750\,167}{26\,460}+\frac
{6583}{2079}+\frac{943}{4004}+\frac{1129}{90\,090}+\frac{217}{437\,580}%
+\frac{2}{138\,567}+\frac{1}{3879\,876}\smallskip
\end{array}
\right)
\]
In the finite sequence of terms that sum to give $\left\vert B_{2n}\right\vert
$ for $n>2$, as displayed above, obviously the second number in the sequence
is just $3/7$ times the first. \ Less obviously, each term in the sequence for
$\left\vert B_{2n}\right\vert $ is greater than the subtotal of all the
smaller terms in that same sequence.

The general result for $\left\vert B_{2n}\right\vert $ is obtained by writing
(\ref{RID2}) as an infinite matrix equation, $\boldsymbol{1}%
=\boldsymbol{M\cdot B}$, where $\boldsymbol{B}$\ is an infinite column of the
even index Bernoulli numbers, $\boldsymbol{1}$\ is an infinite column of $1$s,
and $\boldsymbol{M}_{m,n}=2\left(  -1\right)  ^{m+1}\binom{2n-1}{m}%
\binom{2m+1}{2n}$. \ Computing the inverse for the triangular matrix
$\boldsymbol{M}$ then gives $\boldsymbol{B=M}^{-1}\boldsymbol{\cdot1}$. \ All
terms in a given row of the triangular matrix $\boldsymbol{M}^{-1}$ have the
same sign. \ The ordered terms in the sums above are just the unsigned entries
in $\boldsymbol{M}^{-1}$\ for the $n$th row. \ For example, the first six rows
and columns of the triangular matrices $\boldsymbol{M}$ and $\boldsymbol{M}%
^{-1}$ are given by%
\[
\boldsymbol{M=}\left(
\begin{array}
[c]{cccccc}%
6 & 0 & 0 & 0 & 0 & 0\\
0 & -30 & 0 & 0 & 0 & 0\\
0 & 70 & 140 & 0 & 0 & 0\\
0 & 0 & -840 & -630 & 0 & 0\\
0 & 0 & 924 & 6930 & 2772 & 0\\
0 & 0 & 0 & -18\,018 & -48\,048 & -12\,012
\end{array}
\right)  ,\ \ \boldsymbol{M}^{-1}=\left(
\begin{array}
[c]{cccccc}%
\frac{1}{6}\smallskip & 0 & 0 & 0 & 0 & 0\\
0 & -\frac{1}{30}\smallskip & 0 & 0 & 0 & 0\\
0 & \frac{1}{60}\smallskip & \frac{1}{140} & 0 & 0 & 0\\
0 & -\frac{1}{45}\smallskip & -\frac{1}{105} & -\frac{1}{630} & 0 & 0\\
0 & \frac{1}{20}\smallskip & \frac{3}{140} & \frac{1}{252} & \frac{1}{2772} &
0\\
0 & -\frac{1}{6}\smallskip & -\frac{1}{14} & -\frac{17}{1260} & -\frac{1}{693}
& -\frac{1}{12\,012}%
\end{array}
\right)
\]
hence the previous partitions of $\left\vert B_{2n}\right\vert $ for $n=1$ to
$6$. \ Various things can be said in general about the entries in
$\boldsymbol{M}^{-1}$, e.g., the diagonal is $\left\vert \boldsymbol{M}%
_{n,n}^{-1}\right\vert =\frac{\left(  n!\right)  ^{2}}{\left(  2n+1\right)
!}$, the first sub-diagonal is $\left\vert \boldsymbol{M}_{n\geq2,n-1}%
^{-1}\right\vert =\frac{1}{6}\left(  n-2\right)  \frac{n!\left(  n-1\right)
!}{\left(  2n-1\right)  !}$, the second sub-diagonal is $\left\vert
\boldsymbol{M}_{n\geq3,n-2}^{-1}\right\vert =\frac{7}{360}\left(  n-\frac
{8}{7}\right)  \left(  n-3\right)  \frac{n!\left(  n-2\right)  !}{\left(
2n-3\right)  !}$, etc. \ 

Since the first columns of $\boldsymbol{M}$ and $\boldsymbol{M}^{-1}$ are null
except for the top entry, such that $\left\vert B_{2}\right\vert =1/6$ and
$\left\vert B_{4}\right\vert =1/30$ are both without any partitioning, it
makes sense to begin the set of partitions for $\left\vert B_{2m}\right\vert $
with those for $m=2$. \ One way to express the general result is then
\begin{equation}
\left\vert B_{2m}\right\vert =\sum_{n=2}^{m}b_{m}\left(  n\right)
\ ,\ \ \ m\geq2
\end{equation}
where $b_{m}\left(  n\right)  $ is the non-zero $n^{th}$ \emph{unsigned}
column entry in the $m^{th}$ row of the aforementioned triangular matrix
$\boldsymbol{M}^{-1}$. \ For example, $b_{2}\left(  2\right)  =1/30$,
$b_{3}\left(  2\right)  =1/60$, and $b_{3}\left(  3\right)  =1/140$, with no
other non-zero entries in the 2nd and 3rd rows of $\boldsymbol{M}^{-1}$. \ So
each $\left\vert B_{2m}\right\vert $ is given as a sum of $m-1$ monotonically
falling, positive rational numbers as tabulated above for $2\leq m\leq10$. \ 

For larger values of $m$ and $n$, as presented previously in \cite{TLCSigma},
\begin{equation}
b_{m}\left(  n\right)  =\frac{m!n!}{\left(  2n+1\right)  !}\left(  n-1\right)
q_{m-1-n}\left(  m\right)  \ ,\ \ \ m\geq2
\end{equation}
where again $m$ is the \emph{row} index and $n$ is the \emph{column} index of
$\boldsymbol{M}^{-1}$. \ Here $q_{-1}\left(  n\right)  \equiv\frac{1}{n-1}$
while the $q_{l}\left(  n\right)  $, with $l\geq0$ are $l$th order polynomials
in powers of $n$. \ These polynomials may be obtained sequentially
\cite{TLCSigma} by iteration of
\begin{equation}
q_{l}\left(  n\right)  =\frac{\left(  -1\right)  ^{l}\left(  n-l-3\right)
!}{\left(  2l+3\right)  !\left(  n-2l-3\right)  !}+\sum_{j=0}^{l-1}%
\frac{\left(  -1\right)  ^{l+j+1}\left(  n+j-l-1\right)  !}{\left(
2l+1-2j\right)  !\left(  n+2j-2l-1\right)  !}~q_{j}\left(  n+j-l\right)
\label{recursed}%
\end{equation}
A few explicit examples are \cite{Knuth} \
\begin{equation}
q_{-1}\left(  n\right)  =\frac{1}{n-1}\ ,\ \ \ q_{0}\left(  n\right)
=\frac{1}{6}\ ,\ \ \ q_{1}\left(  n\right)  =\frac{7}{360}~n-\frac{1}%
{45}\ ,\ \ \ q_{2}\left(  n\right)  =\frac{31}{15\,120}~n^{2}-\frac
{89}{15\,120}~n+\frac{1}{315}\ ,\ \ \ \text{etc.}%
\end{equation}
For Bernoulli partitions the $q_{l}\left(  n\right)  $ polynomials are only
needed for $n\geq l+3$. \ 

Now it so happens, for any fixed column $n$, numerical calculations
\cite{TLCunpub} reveal that $b_{m}\left(  n\right)  $ grows without bound as
$m$ increases with asymptotic behavior proportional to that of $\left\vert
B_{2m}\right\vert \sim\frac{2\left(  2m\right)  !}{\left(  2\pi\right)  ^{2m}%
}$. \ Moreover, as the row index $m$ increases, there is a monotonic approach
to an asymptotic ratio
\begin{equation}
\frac{b_{m}\left(  n\right)  }{\left\vert B_{2m}\right\vert }%
\underset{m\rightarrow\infty}{\sim}a\left(  n\right)
\end{equation}
where $a\left(  n\right)  $ is a monotonically falling function of the column
index for $n\geq2$. \ For example, numerically to 10 digits: $a\left(
2\right)  =0.657\,973\,626\,7$, $a\left(  3\right)  =0.281\,988\,697\,2$,
$a\left(  4\right)  =0.0535\,819\,841\,1$, $a\left(  5\right)
=0.005\,985\,334\,466$, etc. \ When augmented by
\href{https://oeis.org/A182448}{https://oeis.org/A182448}, these numerical
revelations lead to the following.

\begin{conjecture}
The asymptotic ratios are $a\left(  n\right)  =p_{n}\left(  \pi^{2}\right)  $
for $n\geq2$ where $p_{n}\left(  x\right)  $ are polynomials:
\begin{equation}
p_{n}\left(  x\right)  =\frac{1}{4n^{2}-1}\sum_{k=0}^{\left\lfloor
n/2\right\rfloor -1}\frac{\left(  -4\right)  ^{k}x^{k+1}}{\left(  2k+1\right)
!}\frac{\Gamma\left(  n-1\right)  \Gamma\left(  2n-2-2k\right)  }%
{\Gamma\left(  n-1-2k\right)  \Gamma\left(  2n-2\right)  }\label{punchline}%
\end{equation}
with $\left\lfloor \ldots\right\rfloor $ representing the floor
function.\footnote{For $n\geq2$, note that $\sum_{k=0}^{\left\lfloor
n/2\right\rfloor -1}\frac{\left(  -4\right)  ^{k}x^{k+1}}{\left(  2k+1\right)
!}\frac{\Gamma\left(  n-1\right)  \Gamma\left(  2n-2-2k\right)  }%
{\Gamma\left(  n-1-2k\right)  \Gamma\left(  2n-2\right)  }=\sum_{k=0}%
^{\left\lfloor \left(  n-1\right)  /2\right\rfloor }\frac{\left(  -2\right)
^{k}x^{k+1}}{\left(  2k+1\right)  !}\left(
%TCIMACRO{\dprod _{l=1}^{k}}%
%BeginExpansion
{\displaystyle\prod_{l=1}^{k}}
%EndExpansion
\frac{\left(  n-k-1-l\right)  }{\left(  2n-1-2l\right)  }\right)  $
\par
as used in earlier versions of this paper, with the empty product equal to $1$
by convention.}
\end{conjecture}

\noindent For example, numerically to 16 digits: $\ p_{2}\left(  \pi
^{2}\right)  =\frac{1}{15}\pi^{2}=0.657\,973\,626\,739\,290\,6\cdots$,
$p_{3}\left(  \pi^{2}\right)  =\frac{1}{35}\pi^{2}%
=0.281\,988\,697\,173\,981\,7\cdots$, $p_{4}\left(  \pi^{2}\right)  =\frac
{1}{63}\pi^{2}\left(  1-\frac{1}{15}\pi^{2}\right)
=0.0535\,819\,841\,082\,941\,2\cdots$, and $p_{5}\left(  \pi^{2}\right)
=\frac{1}{99}\pi^{2}\left(  1-\frac{2}{21}\pi^{2}\right)
=0.005\,985\,334\,466\,027\,733\cdots$, in complete agreement with the
previous numerical results. \ 

\begin{remark}
Here are suggestions for a possible proof of the conjecture. \ In the large
$m$ limit, keeping only leading terms, recast the iteration (\ref{recursed})
as a homogeneous second-order difference equation in $n$ for $a_{n}\equiv
p_{n}\left(  x\right)  $, to obtain%
\begin{equation}
\left(  2n+1\right)  p_{n}\left(  x\right)  =\left(  2n-3\right)
p_{n-1}\left(  x\right)  -\frac{x}{\left(  2n-1\right)  }~p_{n-2}\left(
x\right)  \label{2nd}%
\end{equation}
where $x$ is a parameter related to the normalization of the solution. \ With
chosen initial conditions $p_{1}\left(  x\right)  =0$ and $p_{2}\left(
x\right)  =x/15$, the solution of (\ref{2nd}) is given by (\ref{punchline}).
\ With $p_{1}=0$ it follows from (\ref{2nd}) \ that $7p_{3}\left(  x\right)
=3p_{2}\left(  x\right)  $, which is the same as the exact relation previously
observed in the table of partitions of $\left\vert B_{2n}\right\vert $ for
$3\leq n\leq10$. \ It remains to determine the normalization of $p_{2}$, i.e.,
the value of $x$. \ Allowed values of $x$ must now satisfy two requirements:
\ (1)\ $0\leq p_{n}\left(  x\right)  \leq1$ for all $n\geq2$, and
(2)$\ \sum_{n=2}^{\infty}p_{n}\left(  x\right)  =1$. \ That is to say, the
$p_{n}$ must provide an infinite partition of $1$, as discussed in more detail
below. \ The claim is that only $x=\pi^{2}$ satisfies these
requirements.\medskip
\end{remark}

As re-expressed by Maple and Mathematica, the above polynomials can be
evaluated as generalized hypergeometric functions, namely,%
\begin{equation}
p_{n}\left(  x\right)  =\left(  \frac{x}{4n^{2}-1}\right)  \left.  _{2}%
F_{3}\right.  \left(  \frac{2-n}{2},~\frac{3-n}{2}~;~\frac{3}{2},~\frac
{3-2n}{2},~2-n~;~-x\right)
\end{equation}
Moreover, as suggested by the series form in (\ref{punchline}), the leading
asymptotic behavior for large $n$ then leads back to the sinc function.
\begin{equation}
\left.  _{2}F_{3}\right.  \left(  \frac{2-n}{2},~\frac{3-n}{2}~;~\frac{3}%
{2},~\frac{3-2n}{2},~2-n~;~-x\right)  \underset{n\rightarrow\infty}{\sim}%
\frac{\sin\left(  \sqrt{x}\right)  }{\sqrt{x}}%
\end{equation}
That is somewhat interesting, but more importantly for the partitions of
$\left\vert B_{2m}\right\vert $, even for finite $n$ these particular
hypergeometric functions can be expressed in terms of Bessel functions
\cite{Wolfram} which leads to a more succinct but less transparent closed form
for the conjectured asymptotic ratios.
\begin{equation}
a\left(  n\right)  =p_{n}\left(  \pi^{2}\right)  =\frac{\pi}{\sqrt{2}}%
\frac{n!}{\left(  2n+1\right)  !}\left(  2\pi\right)  ^{n}J_{n-1/2}\left(
\pi\right)  \label{Bessels}%
\end{equation}
Taking (\ref{well-known}) into account then gives
\begin{equation}
b_{m}\left(  n\right)  \underset{m\rightarrow\infty}{\sim}\frac{\left(
2m\right)  !~n!}{\left(  2n+1\right)  !}\left(  2\pi\right)  ^{n-2m}\sqrt
{2}\pi J_{n-1/2}\left(  \pi\right)  \label{bAsymptotic}%
\end{equation}
The right-hand-side of (\ref{bAsymptotic}) is the \textquotedblleft asymptotic
approximant\textquotedblright\ to $b_{m}\left(  n\right)  $ for fixed $n$.

It is perhaps also interesting, even for moderate values of $m$, to compare
numerically\ the asymptotic approximant and exact results. \ For example, with
the exact results as tabulated above,

\begin{center}
$%
\begin{array}
[c]{lll}%
m,\ n\smallskip & \fbox{Exact}\smallskip & \fbox{Asymptotic Approximant}%
\smallskip\\
10,\ 2 & \frac{43\,867}{126}=348.\,151\smallskip & \sim\ 348.\,149\\
10,\ 3 & \frac{43\,867}{294}=149.\,207\smallskip & \sim\ 149.\,207\\
10,\ 4 & \frac{750\,167}{26\,460}=28.\,351\,0\smallskip & \sim\ 28.\,351\,5\\
10,\ 5 & \frac{6583}{2079}=3.\,166\,43\smallskip & \sim\ 3.\,166\,98\\
10,\ 6 & \frac{943}{4004}=0.235\,514\,\smallskip & \sim\ 0.235\,757\,\\
10,\ 7 & \frac{1129}{90\,090}=1.\,253\,19\times10^{-2}\smallskip &
\sim\ 1.\,259\,65\times10^{-2}\\
10,\ 8 & \frac{217}{437\,580}=4.\,959\,09\times10^{-4}\smallskip &
\sim\ 5.\,078\,31\times10^{-4}\\
10,\ 9 & \frac{2}{138\,567}=1.\,443\,35\times10^{-5}\smallskip &
\sim\ 1.\,601\,95\times10^{-5}\\
10,\ 10 & \frac{1}{3879\,876}=2.\,577\,4\times10^{-7} & \sim
\ 4.\,065\,25\times10^{-7}%
\end{array}
$
\end{center}

\noindent As should be expected, the relative error between the approximant
and the exact result increases as the diagonal of $\boldsymbol{M}^{-1}$ is
approached.\smallskip

For a fixed row number $m\geq2$, dividing the non-zero column entries in the
matrix $\boldsymbol{M}^{-1}$\ by the Bernoulli number corresponding to $m$
gives an exact partition of $1$ as a finite sum of $m-1$ monotonic, same-sign,
rational numbers: $\ 1=\sum_{n=2}^{m}\frac{b_{m}\left(  n\right)  }{\left\vert
B_{2m}\right\vert }$.\ \ By extension, using (\ref{Bessels}), the ratios
$a\left(  n\right)  $ also sum to give unity \cite{BesselSum}: \ $1=\sum
_{n=2}^{\infty}a\left(  n\right)  $, where the sum is now infinite because\ in
the limit of infinite row number there are an infinite number of non-zero
entries in the columns of $\boldsymbol{M}^{-1}$. \ In fact the sum over $n$
converges rapidly to $1$ when computed numerically. \ For example,
$1-\sum_{n=2}^{18}p_{n}\left(  \pi^{2}\right)  =2.\,6\times10^{-27}$.

Of course there are an infinite (albeit countable) number of other ways to
partition $1$ into a sum of positive rationals, each less than one and
arranged in a monotonic sequence. \ \emph{But} most such partitions have
\emph{nothing} to do with Bernoulli numbers. \ In contrast to such generic
partitions, in the author's opinion, the scale symmetry origins of the
partitions produced by $\boldsymbol{M}^{-1}$\ make those special partitions
intriguing and worth pursuing.\newpage

In particular, the monotonic finite series representation of $B_{2n}=\left(
-1\right)  ^{n+1}\left\vert B_{2n}\right\vert $ clearly gives a series of
progressively better bounds on $B_{2n}$. \ Such a series of constraints on
$B_{2n}$ might be useful to establish bounds on functions defined as infinite
series whose coefficients involve the Bernoulli numbers \cite{Riesz}. \ That
remains to be shown, but some form of scale symmetry has long been suspected
to be lurking in the Riemann hypothesis \cite{Berry}. \ If so, the progression
of bounds on $B_{2n}$ that follow from the partitions discussed here may
provide some insight.

\textbf{Added Note:} \ More recent work \cite{CV} shows the \emph{exact} ratio
$b_{m}\left(  n\right)  /\left\vert B_{2m}\right\vert $ for any $m\geq n\geq2$
is given by (\ref{punchline}) after the replacement $x^{k}\rightarrow\pi
^{2k}\zeta\left(  2m-2k\right)  /\zeta\left(  2m\right)  $. \ Asymptotic
behavior as $m\rightarrow\infty$ for fixed $k$ is then transparent.

\textbf{Acknowledgements} \ I thank\ C. Vignat for stimulating comments \&
references, and C. Bender \& T.S. Van Kortryk for discussions. \ I received
financial support from the United States Social Security Administration.

\end{document}